\documentclass[12pt]{article}
\usepackage[pctex32]{graphics}
\usepackage{amsthm,amsmath,amssymb,amscd,verbatim,epsfig}
\usepackage{graphics,pifont}
\usepackage{amssymb}          
\usepackage{graphicx}
\newcommand{\beq}{\begin{equation}}
\newcommand{\eeq}{\end{equation}}

\date{}

\newcommand{\f}{\frac}
\newcommand{\ra}{\rightarrow}

\newcommand{\supp}{\supset}

\newcommand{\sq}{$\blacksquare$}

\begin{document}

\title{A\ Simple\ Proof\ of\ Sharkovsky's\ Theorem}
\author{Bau-Sen Du \\ [.3cm]
Institute of Mathematics \\
Academia Sinica \\
Taipei 11529, Taiwan \\
dubs@math.sinica.edu.tw \\}
\maketitle

\begin{abstract}
In this note, we present a simple directed graph proof of Sharkovsky's theorem.  
\end{abstract}


\section{Introduction.}
Throughout this note, $I$ is a compact interval, and $f : I \ra I$ is a continuous map.   For each integer $n \ge 1$, let $f^n$ be defined by: $f^1 = f$ and $f^n = f \circ f^{n-1}$ when $n \ge 2$.  For $y$ in $I$, we call the set $O_f(y)=\{ \, f^k(y) \, \big| \, k \ge 0 \, \}$ the {\it{orbit}} of $y$ (under $f$) and call $y$ a {\it{periodic point}} of $f$ with least period $m$ (or a period-$m$ point of $f$) if $f^m(y) = y$ and $f^i(y) \ne y$ when $0 < i < m$.  If $f(y) = y$, then we call $y$ a {\it{fixed point}} of $f$.  It is clear that every $f$ of the type in question has fixed points.

For discrete dynamical systems defined by iterated interval maps, one of the most remarkable results is Sharkovsky's theorem {\bf{\cite{mi}}}, {\bf{\cite{sh}}}.  It states that, if $f$ has a period-$m$ point, then $f$ also has a period-$n$ point precisely when $m \prec n$ in the following Sharkovsky's ordering : $$3 \prec 5 \prec 7 \prec \cdots \prec 2 \cdot 3 \prec 2 \cdot 5 \prec 2 \cdot 7 \prec \cdots \prec 2^2 \cdot 3 \prec 2^2 \cdot 5 \prec 2^2 \cdot 7 \prec \cdots \prec 2^3 \prec 2^2 \prec 2 \prec 1.$$  It is well known (see {\bf{\cite{str}}}) that the sufficiency of Sharkovsky's theorem can be derived from the following three statements: (a) if $f$ has a periodic point of least period greater than 2, then $f$ also has a periodic point of least period 2; (b) if $f$ has a periodic point of odd period $m \ge 3$, then $f$ also has a periodic point of least period $n$ for every integer $n$ such that $n \ge m+1$; (c) if $f$ has a periodic point of odd period $m \ge 3$, then $f$ also has periodic points of all even periods.  The difficulty of proving the sufficiency of Sharkovsky's theorem lies in proving (c), where most proofs involve the structures of the so-called \v Stefan cycles {\bf{\cite{al}}}, {\bf{\cite{bl}}}, {\bf{\cite{ste}}}.  In this note, we give a unified proof of (b) and (c) that does not involve \v Stefan cycles.  We also give a different proof of (a) {\bf{\cite{bb}}}, {\bf{\cite{str}}}.  For the sake of completeness, we end the paper with a proof of Sharkovsky's theorem.

In proving Sharkovsky's theorem, we need the following result [{\bf 3}, p.12].
\medskip

\noindent
{\bf Lemma 1.}
{\it Let $k, m, n$, and $s$ be positive integers.  Then the following statements hold:
\begin{itemize}
\item[\rm{(1)}]
If $y$ is a periodic point of $f$ with least period $m$, then it is a periodic point of $f^n$ with least period $m/(m, n)$, where $(m, n)$ is the greatest common divisor of $m$ and $n$. 

\item[\rm{(2)}]
If $y$ is a periodic point of $f^n$ with least period $k$, then it is a periodic point of $f$ with least period $kn/s$, where $s$ divides $n$ and is relatively prime to $k$. 
\end{itemize}}

\noindent
\section{A proof of (a).}
We need the following result, which can also be used to show that $f$ has no period-2 point if and only if for each point $c$ of $I$ the iterates $f^n(c)$ converge to a fixed point of $f$ [{\bf 3}, p.121].

\noindent
{\bf Lemma 2.}
{\it If $c$ and $d$ are points of $I$ such that $f(d) \le c < d \le f(c)$, then $f$ has a periodic point of least period 2. } 

\noindent
{\it Proof.}
Write $I = [a, b]$.  Let $w = \min \{ c \le x \le d \, \big| \, f(x) = x \, \}$, and let $v$ be a point in $[c, w]$ with $f(v) = d$.  Then, $f^2(v) = f(d) \le c \le v$.  If $f$ has no fixed point in $[a, c]$, then it fixes no point of $[a, v]$.  Since $f^2(a) \ge a$, it follows that $f$ has a periodic point with least period 2 in $[a, v]$.  If $f$ has a fixed point in $[a, c]$, let $t = \max \{ a \le x < c \, \big| \, f(x) = x \, \}$.  Then $f$ has no fixed point in $(t, v]$.  Let $u$ be a point in $[t, c]$ with $f(u) = c$.  Then $f^2(u) = f(c) \ge d > u$.  Since $f^2(v) \le v$, we infer that $f^2(y) = y$ for some $y$ in $[u, v]$.  Because $f$ has no fixed point in $[u, v]$, $y$ is a periodic point of $f$ with least period 2.
\hfill\sq

\noindent
{\bf Proposition 3.}
{\it If $f$ has a periodic point of least period $m$ larger than 2, then $f$ also has a periodic point of least period 2.}

\noindent
{\it Proof.}
Let $P = \{ \, x_i \, \big| \, 1 \le i \le m \, \}$, with $x_1 < x_2 < \cdots < x_m$, be a period-$m$ orbit of $f$.  Since $x_1 < f(x_1)$ and $f(x_m) < x_m$, there exists an integer $s$ satisfying $1 \le s \le m-1$ such that $x_s = \max \{ \, x \in P \, \big| \, x < f(x) \, \}$.  It is clear that $x_{s+1} \le f(x_s)$ and $f(x_{s+1}) \le x_s$.  By Lemma 2, $f$ has a periodic point of least period 2.
\hfill\sq

\noindent
\section{A unified proof of (b) and (c).}
If there are closed subintervals $J_0$, $J_1$, $\cdots$, $J_{n-1}, J_n$ of $I$ with $J_n = J_0$ such that $f(J_i) \supp J_{i+1}$ for $i = 0, 1, \cdots, n-1$, then we say that $J_0J_1 \cdots J_{n-1}J_0$ is a {\it{cycle of length}} $n$.  We require the following result. 

\noindent
{\bf Lemma 4.}
{\it If $J_0J_1J_2 \cdots J_{n-1}J_0$ is a cycle of length $n$, then there exists a periodic point $y$ of $f$ such that $f^i(y)$ belongs to $J_i$ for $i = 0, 1, \cdots, n-1$ and $f^n(y) = y$.}  

We now give a simple unified proof of (b) and (c).

\noindent
{\bf Proposition 5.}
{\it If $f$ has a periodic point of least period $m$ with $m \ge 3$ and odd, then $f$ has periodic points of all even periods.  Furthermore, $f$ has a periodic point of least period $n$ for each integer $n$ with $n \ge m+1$.}

\noindent
{\it Proof.}
Let $P = \{ \, x_i \, \big| \, 1 \le i \le m \, \}$, with $x_1 < x_2 < \cdots < x_m$, be a period-$m$ orbit of $f$.  Let $x_s = \max \{ \, x \in P \, \big| \, x < f(x) \, \}$.  Then $x_{s+1} \le f(x_s)$ and $f(x_{s+1}) \le x_s$, so $f$ has a fixed point $z$ in $[x_s, x_{s+1}]$.  Since $m$ is odd, for some integer $t$ such that $1 \le t \le m-1$ and $t \ne s$ the points $f(x_t)$ and $f(x_{t+1})$ lie on opposite sides of $z$.  Thus $f([x_t, x_{t+1}])$ $\supp [x_s, x_{s+1}]$.  For simplicity, we assume that $x_t < x_s$.  If $x_{s+1} \le x_t$, the proof is similar.  Let $q$ be the smallest positive integer such that $f^q(x_s) \le x_t$.  Then $2 \le q \le m-1$.

First assume that $m = 3$.  Without loss of generality, we assume that $f(x_1)=x_2$, $f(x_2)=x_3$, and $f(x_3)=x_1$.  Let $J_0=[x_1,x_2]$ and $J_1=[x_2,x_3]$.  For any $n \ge 2$, we can apply Lemma 4 to the cycle $J_0J_1J_1 \cdots J_1J_0$ of length $n$ to obtain a period-$n$ point.  Accordingly, if $f$ has a period-3 point, then $f$ has periodic points of all periods.  Now assume that $m > 3$.  Since $q$ is the smallest positive integer such that $f^q(x_s) \le x_t$, $x_{t+1} \le f^i(x_s)$ whenever $1 \le i \le q-1$.  If $x_{t+1}$ $\le f^{q-1}(x_s)$ $< x_s$, Lemma 4 applies to the cycle $$[x_t, f^{q-1}(x_s)][f^{q-1}(x_s), z][f^{q-1}(x_s), z][x_t, f^{q-1}(x_s)]$$ and establishes the existence of a period-3 point of $f$.  If $f^{q-1}(x_s) = x_{s+1}$, we can apply Lemma 4 to the cycle $$[z, x_{s+1}][x_t, x_{t+1}][x_s, x_{s+1}][z, x_{s+1}]$$ to obtain a period-3 point of $f$.  

We proceed assuming that $x_{s+1} < f^{q-1}(x_s)$.  If $k = \min \{$ \, $1 \le i \le q-1$ $\, \big| \, f^{q-1}(x_s) \le f^i(x_s) \, \}$, then $x_{t+1} \le f^{k-1}(x_s) < f^{q-1}(x_s)$, so either $x_{s+1}$ $\le f^{k-1}(x_s)$ $< f^{q-1}(x_s)$ or $x_{t+1} \le f^{k-1}(x_s) \le x_s$.  If $x_{s+1} \le f^{k-1}(x_s) < f^{q-1}(x_s)$, we can invoke Lemma 4 for the cycle $$[f^{k-1}(x_s), f^{q-1}(x_s)][z, f^{k-1}(x_s)][z, f^{k-1}(x_s)][f^{k-1}(x_s), f^{q-1}(x_s)]$$ to obtain a period-3 point of $f$.  If $x_{t+1}$ $\le f^{k-1}(x_s)$ $\le x_s$ $(< z$ $< f^{q-1}(x_s))$, we choose $u$ in $[x_t, x_{t+1}]$ such that $f(u) = z$, pick $w$ in $[z, f^{q-1}(x_s)]$ with $f(w) = u$, and let $v$ in $[f^{k-1}(x_s), z]$ be a point such that $f(v) = w$.  By applying Lemma 4 to the cycle $[u, v][z, w][u,v]$ and, for every {\it{even}} integer $n \ge 4$, to the cycle $$[u, v]([z, w][v, z])^{\f {n-2}2}[z, w][u, v]$$ (here $([z, w][v, z])^{\f {n-2}2}$ represents $(n-2)/2$ copies of $[z, w][v, z]$) of length $n$, we conclude that $f$ has periodic points of all even periods.  On the other hand, let $J_i = [z : f^i(x_s)]$ for $i = 0, 1, \cdots, q-1$, where $[a : b]$ denotes the closed interval with $a$ and $b$ as endpoints.  For any $n \ge m+1$, we appeal to Lemma 4 to the cycle of length $n$ \, $J_0J_1$ $\cdots$ $J_{k-1}J_{q-1}[x_t, x_{t+1}]J \cdots JJ_0$, where $J = [x_s, x_{s+1}]$, to confirm the existence of a period-$n$ point.
\hfill\sq

\noindent
\section{A proof of Sharkovsky's theorem.}  We now combine (a), (b), (c), and Lemma 1 to prove Sharkovsky's theorem.  

\noindent
{\bf Theorem 6 (Sharkovsky).}
{\it Assume that $f : I \ra I$ is a continuous map.  If $f$ has a period-$m$ point, then $f$ also has a period-$n$ point precisely when $m \prec n$ in the Sharkovsky's ordering defined as in Section 1.}

\noindent
{\it Proof.}
By (b) and (c), we have $3 \prec 5 \prec 7 \prec \cdots \prec 2 \cdot 3$.  If $f$ has period-$(2 \cdot m)$ points with $m \ge 3$ and odd, then $f^2$ has period-$m$ points.  By (b), $f^2$ has period-$(m+2)$ points, which by Lemma 1(2) implies that $f$ has either period-$(m+2)$ points or period-$(2 \cdot (m + 2))$ points.  If $f$ has period-$(m + 2)$ points, then by (b) $f$ also has period-$(2 \cdot (m+2))$ points.  In either case, $f$ has period-$(2 \cdot (m+2))$ points.  On the other hand, by (c) $f^2$ has period-6 points and hence, by Lemma 1(2), $f$ has period-$(2^2 \cdot 3)$ points.  Now if $f$ has period-$(2^k \cdot m)$ points with $m \ge 3$ and odd and if $k \ge 2$, then by Lemma 1(1) $f^{2^{k-1}}$ has period-$(2 \cdot m)$ points.  It follows from what we have just proved that $f^{2^{k-1}}$ has period-$(2 \cdot (m+2))$ points and period-$(2^2 \cdot 3)$ points.  In view of Lemma 1(2), $f$ has period-$(2^k \cdot (m+2))$ points and period-$(2^{k+1} \cdot 3)$ points.  Furthermore, because $f$ has period-$(2^k \cdot m)$ points, $f^{2^k}$ has period-$m$ points.  By (b), $f^{2^k}$ has period-$2^n$ points as long as $2^n > m$, so by Lemma 1(2) $f$ has period-$(2^{k+n})$ points for all integers $n$ such that $2^n > m$.  Finally, if $f$ has period-$2^i$ points for some integer $i \ge 2$, then $f^{2^{i-2}}$ has period-4 points.  As a result of (a), $f^{2^{i-2}}$ has period-2 points, ensuring that $f$ has period-$2^{i-1}$ points.  This proves the sufficiency of Sharkovsky's theorem. 

For the converse, it suffices to assume that $I = [0, 1]$.  Let $T(x) = 1 - |2x - 1|$ be the tent map on $I$.  Then for any $k \ge 1$ the equation $T^k(x) = x$ has exactly $2^k$ distinct solutions in $I$.  It follows that $T$ has finitely many period-$k$ orbits.  Among these period-$k$ orbits, let $P_k$ be one with the smallest diameter $\max P_k - \min P_k$.  For any $x$ in $I$, let $T_k(x) = \min P_k$ if $T(x) \le \min P_k$, $T_k(x) = \max P_k$ if $T(x) \ge \max P_k$, and $T_k(x) = T(x)$ if $\min P_k \le T(x) \le \max P_k$.  It is then easy to see that $T_k$ has exactly one period-$k$ orbit, i.e., $P_k$, and no period-$j$ orbit for any $j$ with $j \prec k$ in the Sharkovsky's ordering (see also {\bf{\cite{al}}}, pp. 32-34]).  Now let $Q_3$ be any period-3 orbit of $T$ of minimal diameter.  Then $[\min Q_3, \max Q_3]$ contains finitely many period-6 orbits of $T$.  If $Q_6$ is one of smallest diameter, then $[\min Q_6, \max Q_6]$ contains finitely many period-12 orbits of $T$.  We choose one, say $Q_{12}$, of minimal diameter and continue the process inductively.  Let $q_0 = \sup \{\min Q_{2^n \cdot 3} \, \big| \, n \ge 0 \}$ and $q_1 = \inf \{ \max Q_{2^n \cdot 3} \, \big| \, n \ge 0 \}$.  Let $T_{\infty}(x) = q_0$ if $T(x) \le q_0$, $T_{\infty}(x) = q_1$ if $T(x) \ge q_1$, and $T_{\infty}(x) = T(x)$ if $q_0 \le T(x) \le q_1$.  Then it is easy to check that $T_{\infty}$ has periodic points of least period $2^n$ for each $n \ge 0$, but has no periodic points of any other periods.  This establishes the other direction in Sharkovsky's theorem.
\hfill\sq

\noindent
{\bf Remark.} Our method can also be used to prove that if $f$ has a periodic point of odd period $m > 1$, but no periodic points of odd period strictly between 1 and $m$ then any periodic orbit of odd period $m$ must be a \v Stefan orbit (cf. {\bf{\cite{bu}}}).

\noindent
{\bf Acknowledgments.}  I would like to thank M. Misiurewicz, A. N. Sharkovsky, and the referee for many constructive suggestions that led to improvements in this note.

\end{document}